\documentclass[11pt]{amsart}
\usepackage{amssymb}
\usepackage{graphicx}
\usepackage{slashed}
\usepackage[colorlinks]{hyperref} 
\usepackage[all,cmtip]{xy}

\newtheorem{thm}{Theorem}[subsection]

\newtheorem{thm-def}[thm]{Theorem-Definition}
\newtheorem{lem}[thm]{Lemma}

\theoremstyle{remark}
\newtheorem{rem}[thm]{Remark}
\theoremstyle{remark}

\numberwithin{equation}{subsection}

\newcommand{\vacB}{|0\rangle}
\newcommand{\vac}{\vacB}

\newcommand{\BA}{\mathbb A}
\newcommand{\CC}{\mathbb C}
\newcommand{\BC}{\mathbb C}

\newcommand{\ZZ}{\mathbb Z}
\newcommand{\BZ}{\mathbb Z}

\newcommand{\cA}{\mathcal{A}}

\newcommand{\cC}{\mathcal{C}}

\newcommand{\cH}{{\mathcal{H}}}

\newcommand{\cL}{{\mathcal{L}}}

\newcommand{\cR}{{\mathcal{R}}}

\newcommand{\cV}{{\mathcal{V}}}
\newcommand{\cW}{{\mathcal{W}}}

\begin{document}

\title{one example of a chiral lie group}

\author{ A.~Linshaw\and F.~Malikov}

\thanks{A.~Linshaw is partially supported by Simons Foundation Grant \#318755. }

\thanks{F.~Malikov is partially supported by an NSF grant}

\maketitle

\noindent

\begin{abstract}
We quantize the Khesin-Zakharevich Poisson-Lie group of pseudo-differential symbols.
\end{abstract}

\section{introduction and statement of results}
\label{introduction}
\subsection{ }
\label{intro-1}
A Lie group $G$ is called a {\em Poisson-Lie} group if (cf. \cite{LW})

(i) $G$ is a Poisson manifold and,

(ii) the multiplication $m: G\times G\rightarrow G$ is a Poisson morphism, where $G\times G$ is equipped with the product Poisson structure.

If $G$ is an algebraic group then these conditions become

(i') the structure ring $\BC[G]$ is a Poisson algebra, and

(ii') $m^\#: \BC[G]\rightarrow\BC[G]\otimes\BC[G]$ is a Poisson algebra morphism, where the Poisson algebra structure on $\BC[G]\otimes\BC[G]$ is defined by the condition
that $\BC[G]\otimes1$ and $1\otimes\BC[G]$ Poisson commute with each other and are each naturally isomorphic to $\BC[G]$ as Poisson algebras.

B.~Khesin and I.~Zakharevich \cite{KhZ} suggested the following remarkable example of a Poisson-Lie group.  A  {\em pseudo-differential symbol} of order $\lambda\in\BC$ is an expression of the form
\[
\partial^\lambda+U_1\partial^{\lambda-1}+U_2\partial^{\lambda-2}+\cdots+U_{j}\partial^{\lambda-j}+\cdots,
\]
where $U_j\in\BC((z))$, $j\geq 1$. Pseudo-differential symbols can be multiplied if one postulates the relation
\[
\partial^\mu U=U\partial^{\mu}+\mu U^{(1)}\partial^{\mu-1}+{\mu\choose 2}U^{(2)}\partial^{\mu-2}+\cdots=\sum_{j=0}^\infty {\mu\choose j}U^{(j)}\partial^{\mu-j}.
\]
Using this relation repeatedly one obtains
\begin{equation}
\label{def-of-mult}
(\partial^\lambda+\sum_{j=1}^\infty W_j\partial^{\lambda-j})\cdot(\partial^\mu+\sum_{j=1}^\infty U_j\partial^{\mu-j})=\partial^{\lambda+\mu}+\sum_{j=1}^\infty 
P_j\partial^{\lambda+\mu-j},
\end{equation}
where each $P_j$ is a differential polynomial in $W_1, U_1, W_2, U_2,...$ with coefficients in $\BC[\lambda]$. For example,
\[
P_1=W_1+U_1, \; P_2=W_2+W_1 U_1+\lambda U_1^{(1)}+U_2,\text{ etc.}
\]

The geometric series formula shows how to find the multiplicative inverse:
\begin{equation}
\label{def-of-inv}
(\partial^\lambda+\sum_{j=1}^\infty U_j\partial^{\lambda-j})^{-1}=\partial^{-\lambda}\sum_{n=0}^\infty(-\sum_{j=1}^\infty U_j\partial^{-j})^n.
\end{equation}

This shows that the set of all pseudo-differential symbols of all complex orders carries a group  structure.
Denote this group by $G_{\Psi DO}$.  In fact, $G_{\Psi DO}$ is a group object in the category of ind-schemes, but we will not be too concerned with this structure. Suffice it to say
that for each differential polynomial $P$ in $U_{j}$, $j\geq 1$, $\oint z^n P\,dz$ is a function on $G_{\Psi DO}$; in particular, the set of functions $\{\oint z^n U_{j}\,dz\}_{j\geq 1,n\in\ZZ}$ along
with $\lambda$ has the meaning of a coordinate system on $G_{\Psi DO}$ and the whole of $\BC[G_{\Psi DO}]$ is a topological algebra ``generated'' by them. For example, the equation $\lambda=\lambda_0$ defines the set of all order $\lambda_0$ pseudo-differential symbols, and the infinite system of equations $\lambda-n=U_{n+1}=U_{n+2}=\cdots=0$ defines the more familiar set of  ordinary order $n$ differential operators, to be denoted
$\mbox{DO}_n$.

More important for us will be the  algebra of differential polynomials $J_\infty\BC[U_j, j\geq 1][\lambda]$, over $\BC[\lambda]$. If we unburden the notation by letting $\BC[G_{\Psi DO}]^c$ stand for $J_\infty\BC[W_j, j\geq 1][\lambda]$, the above multiplication and inverse maps (\ref{def-of-mult},\ref{def-of-inv}) become
associative algebra morphisms
\begin{equation}
\label{coiss-str-morph}
\BC[G_{\Psi DO}]^c\longrightarrow \BC[G_{\Psi DO}]^c\otimes\BC[G_{\Psi DO}]^c,
\end{equation}
\begin{equation}
\label{coiss-alg-inv-prelim}
\BC[G_{\Psi DO}]^c\longrightarrow \BC[G_{\Psi DO}]^c.
\end{equation}

The Khesin-Zakharevich discovery was that the celebrated Gelfand-Dickey bracket defines on $G_{\Psi DO}$ a Poisson-Lie group structure. More precisely -- and in a  different
terminology -- their result is as follows:

\begin{thm}
\label{kh-zakh-thm}
\begin{sloppypar}
{\em (i) $\BC[G_{\Psi DO}]^c$ carries a coisson  (see \cite{BD}, a.k.a. vertex Poisson \cite{FBZ})  algebra  structure  s.t. (\ref{coiss-str-morph}) is a coisson algebra morphism;

(ii) The differential ideals generated by $\lambda-\lambda_0$, $\lambda_0\in\BC$, and by $\{(\lambda-n, U_{j}, j\geq n+1\}$, $n\in \{1,2,3,...\}$, are coisson. The usual Poisson algebra attached to
the coisson algebra 
\[
\BC[G_{\Psi DO}]^c/(\lambda-n, U_j,j\geq n+1)=:\BC[\mbox{DO}_n]^c
\]
 is isomorphic to the  Gelfand-Dickey Poisson algebra structure on $\mbox{DO}_n$, \cite{GD}.}
\end{sloppypar}

\end{thm}

\begin{rem}
\label{intro-rem-1}
If $n=1$, the coisson algebra  $\BC[\mbox{DO}_1]^c$ is nothing but the Heisenberg algebra on 1 generator:  $\BC[\mbox{DO}_1]=J_\infty\BC[I]$ with coisson bracket
\[
\{I(z),I(w)\}=\delta(z-w)^{(1)}.
\]
As $n$ increases the coisson algebra  $\BC[\mbox{DO}_n]$ becomes more complicated and is best described by the Miura transform, see sect.~\ref{intro-3}  below.
\end{rem}

\begin{rem} 
\label{intro-remark} It follows from the definition \cite{LW} that for any Poisson-Lie group the counit $\BC[G]\rightarrow\BC$ is a Poisson algebra morphism and the inverse $G\rightarrow G$
induces the Poisson algebra isomorphism $\BC[G]\rightarrow\BC[G]^-$, where $\BC[G]^{-}$ means the commutative algebra $\BC[G]$ with Poisson bracket equal to minus
Poisson bracket on $\BC[G]$. The same is true in the vertex Poisson situation here: the counit
\begin{equation}
\label{coiss-counit}
\BC[G_{\Psi DO}]^c\longrightarrow \BC,\; \lambda\mapsto 0, U_j\mapsto 0\mbox{ if } j\geq0.
\end{equation}
is a coisson algebra morphism; similarly,  (\ref{coiss-alg-inv-prelim}) is in fact a coisson algebra isomorphism
\begin{equation}
\label{coiss-alg-inv}
\BC[G_{\Psi DO}]^c\longrightarrow \BC[G_{\Psi DO}]^{c,-},
\end{equation}
where the superscript ``-'' in $\BC^c[G_{\Psi DO}]^{-}$ is easiest to explain in the Beilinson-Drinfeld language \cite{BD}: if the coisson bracket on $\BC^c[G_{\Psi DO}]$ is
 defined as a *-operation $\{.,.\}\in P^*_2(\{\BC^c[G_{\Psi DO}],\BC^c[G_{\Psi DO}]\},\BC^c[G_{\Psi DO}])$, then the bracket on  $\BC^c[G_{\Psi DO}]^{-}$ is $-\{.,.\}$.

It is rather clear that the listed properties of $\BC[G_{\Psi DO}]^c$ earn $G_{\Psi DO}$ the name of {\em coisson-Lie group}. One can then attempt to make a definition
of the latter concept as follows: call a coisson algebra $\cR$ a coisson-Lie group if it carries three coisson algebra morphisms
$\cR\rightarrow \cR\otimes\cR$, $\cR\rightarrow\BC$, and $\cR\rightarrow\cR^-$ s.t. a list of obvious axioms is satisfied. Furthermore, mimicking the discussion in 2.4.9, 2.6.2 of \cite{BD}, one can suggest that the corresponding group object is the ind-scheme of the horizontal sections of $\mbox{Spec}\cR((z))$ over the formal disc $\mbox{Spec}\BC((z))$.
 We will not pursue this here, and restrict ourselves to the following simple remark: if $G$ is a Poisson-Lie
group, then the jet group $J_\infty G$ is coisson-Lie; at least the jet algebra $J_\infty\BC[G]$ is automatically a coisson algebra, and the multiplication map gives a coisson algebra morphism
\[
J_\infty\BC[G]\longrightarrow J_\infty\BC[G]\otimes J_\infty\BC[G].
\]
The example of $\BC[G_{\Psi DO}]$ is more sophisticated.
\end{rem}

\subsection{ }
\label{intro-3}

Our main result consists in the quantization of the Khesin-Zakharevich theorem~\ref{kh-zakh-thm}. Since the category of vertex algebras $\cV ert\cA lg$ carries a natural monoidal structure, one defines
(see \cite{BD}, 3.4.16) a {\em vertex bialgebra } to be a coassociative coalgebra object of $\cV ert\cA lg$, i.e., a $V\in\cV ert\cA lg$ along with  vertex algebra morphisms
$V\rightarrow V\otimes V$ (comultiplication) and $V\rightarrow\BC$ (counit)  s.t. the usual identities hold true. 

The desired quantization will be constructed by a sort of analytic continuation of the classic Lukyanov construction \cite{Luk} of the quantization of the Gelfand-Dickey bracket on the space
of order $n$ differential operators. Lukyanov's construction, also known as the quantum  Miura transform, is as follows. Let $\cH(\hbar)=J_\infty[I_1,...,I_n][\hbar]$; this is an algebra of differential polynomials in the indicated variables over the algebra $\BC[\hbar]$. We give this space a Heisenberg vertex algebra structure by defining the OPEs as follows:
\begin{equation}
\label{def-of-heis}
I_i(z)I_j(w)=\delta_{i,j}\frac{\hbar}{(z-w)^2}+\cdots
\end{equation}
$\cH(\hbar)$, then, is a vertex $\BC[\hbar]$-algebra. More generally, given a unital, commutative, associative algebra $A$ with zero derivation (which makes $A$ a commutative vertex algebra) and a vertex algebra morphism
$A\rightarrow V$, one says that $V$ is a {\em vertex $A$-algebra}. Vertex $A$-algebras are a staple of the present paper.

Consider the formal expression $(\partial_z+I_1(z))(\partial_z+I_2(z))\cdots (\partial_z+I_n(z))$, its purpose being to encode the following expansion
\[
(\partial_z+I_1(z))(\partial_z+I_2(z))\cdots (\partial_z+I_n(z))=\partial_z^n+ ^nU_1(z)\partial_z^{n-1}+\cdots+ ^nU_j(z)\partial_z^{n-j}+\cdots+ ^nU_n(z).
\]

Denote by $\cW_n(\hbar)$ the vertex subalgebra of $\cH(\hbar)$ generated by the fields $^nU_1(z)$,...,$^nU_n(z)$.  Lukyanov proved in \cite{Luk} that this vertex algebra is a quantization
of the Gelfand-Dickey coisson algebra of functions on the space of order $n$ differential operators.  Notice that the fact that $\cW_n(\hbar)/(\hbar)$ is the Gelfand-Dickey algebra is also known
as the {\em Miura transform} and is, from our point of view, simply the statement that the product map
\[
\mbox{DO}_1^{\times n}\longrightarrow\mbox{DO}_n
\]
is Poisson, see Theorem~\ref{kh-zakh-thm} (ii) and Remark~\ref{intro-rem-1}. Lukyanov's result, therefore, is that for any $i,j$ the OPE $^nU_i(z)^nU_j(w)$ involves only normally ordered products of the fields $^nU_\bullet(w)$ and their derivatives,
a highly nontrivial assertion, which has no analogue of comparable simplicity for  most, if not all, other simple Lie algebras.

We conclude this reminder on \cite{Luk} by making the following observation: given arbitrary positive integer $m$ and $n$ splitting the product
\begin{eqnarray}
& &(\partial_z+I_1(z))(\partial_z+I_2(z))\cdots (\partial_z+I_{m+n}(z))=\nonumber\\
& &(\partial_z+I_1(z))(\partial_z+I_2(z))\cdots (\partial_z+I_{m}(z))
(\partial_z+I_{m+1}(z))(\partial_z+I_{m+2}(z))\cdots (\partial_z+I_{m+n}(z))\nonumber
\end{eqnarray}
defines a vertex algebra morphism
\begin{equation}
\label{def-of-alm-comulti}
\cW_{m+n}(\hbar)\longrightarrow\cW_m(\hbar)\otimes\cW_n(\hbar).
\end{equation}
The reader may want to consult \cite{Gen}, Theorem B for a different perspective.

\subsection{ }
\label{intro-4}

We would like to analyze how a vertex algebra structure on $\cW_n(\hbar)$ depends on $n$. Consider the direct system of vector spaces (and by no means vertex algebras!)
\[
\cW_m(\hbar)\rightarrow\cW_n(\hbar),\; ^mU_i\mapsto\, ^nU_i\mbox{ if }m\leq n.
\]

We shall show that for fixed $i,j$ the OPE $^nU_i(z)^nU_j(w)$ is a polynomial in $n$ if $n>N_{ij}$ for some integer $N_{ij}$. If we denote this polynomial by $P_{ij}(n)$, then we can define on  $\lim_{n\rightarrow\infty}\cW_n(\hbar)$ 
a vertex $\BC[\nu,\hbar]$-algebra structure by letting  $U_i(z)U_j(w)\stackrel{\mbox{def}}{=}P_{ij}(\nu)$; here $U_\bullet(z)$ stands for the class of $^nU_\bullet(z)$. Denote thus obtained vertex algebra by $\cL(\nu,\hbar)$. 

The following is the main result of this paper:

\begin{thm}
 \label{main-res}
 (i) The assignment 
 \[ 
 U_j\mapsto\sum_{0\leq a+b\leq i+j}{\nu\otimes1-a\choose j-a-b}U_a\otimes U_b^{(j-a-b)}
 \]
 uniquely extends to a vertex algebra morphism
 \[
 \Delta:\; \cL(\nu,\hbar)\longrightarrow\cL(\nu,\hbar)\otimes\cL(\nu,\hbar).
 \]
 (ii) The set $\{\nu, U_j,\; j\geq 1\}$ generates a codimension 1 vertex ideal of $\cL(\nu,\hbar)$. The composition
 \[
 \BC\longrightarrow \cL(\nu,\hbar)/(\nu, U_j,\; j\geq 1),\; x\mapsto x\vac\mbox{ mod } (\nu, U_j,\; j\geq 1)
 \]
 is an isomorphism; here $\vac\subset\cL(\nu,\hbar)$ is the unit.
 
 (iii) The collection comprising the vertex algebra $\cL(\nu,\hbar)$ along with $\Delta$ from (i) and the composition
 \[
 \eta: \cL(\nu,\hbar)\longrightarrow \cL(\nu,\hbar)/(\nu, U_j,\; j\geq 1)\longrightarrow\BC,
 \]
 where the rightmost arrow is the inverse of the one in (ii), is a vertex bi-algebra. This bi-algebra is a quantization of $\BC[G_{\Psi DO}]^c$ so that the latter is obtained as $\hbar\rightarrow 0$.
 
 (iv) There are vertex algebra surjections 
 \[
 \cL(\nu,\hbar)/(\nu-n)\longrightarrow \cW_n(\hbar),\; n\geq 3,
 \]
 so that the comultiplication $\Delta$, which appeared in item (i),  descends to $\cW_{m+n}(\hbar)\longrightarrow\cW_m(\hbar)\otimes\cW_n(\hbar)$ that was introduced in (\ref{def-of-alm-comulti}).
 \end{thm}
 
 The algebra $\cL(\nu,\hbar)$ appears in \cite{Pro},  conjecturally and inexplicitly, so far as we can see, but with lots of concrete computations including
  the formulas for  comultiplication.

 \subsection{ }
 \label{intro-linshaw-alg}
 Instrumental in the proof of Theorem~\ref{main-res} is another vertex $\BC[\lambda,c]$-algebra, $\cW(\lambda,c)$, introduced by one of us \cite{L} for a similar, yet different, purpose.  
 This algebra solves a certain classification problem and is, in fact, almost a universal object in a certain subcategory of vertex algebras. Namely, for any associative, commutative 
 algebra $A$ a vertex $A$-algebra $V$
 satisfying a number of conditions defines an algebra homomorphism $f:\; \BC[\lambda, c]\rightarrow \tilde{A}$ and a surjection 
 $\tilde{A}\otimes_{\BC[\lambda,c]}\cW(\lambda,c)\rightarrow \tilde{A}\otimes_{A}V$, for some quadratic Galois localization $A\rightarrow\tilde{A}$.
 In our case, this gives us a morphism $ \cW(\lambda,c)\rightarrow \tilde{A}\otimes_{\BC[\nu,\hbar]}\cL(\nu,\hbar)$. Found in \cite{L} is considerable information on what vertex algebras can be obtained as quotients of $\cW(\lambda,c)$. In particular, there are explicitly written down polynomials $p_{n+1}\in\BC[\lambda,c]$ such that the quotient
 $\cW(\lambda,c)/(p_{n+1})$ affords a morphism to the ordinary $W$-algebra, $\cW(sl_n)$. This information suffices to prove Theorem~\ref{main-res}(iv).
 
 In fact, we show that the tensor product $\cH\otimes\cW(\lambda,c)$ and $\cL(\nu,\hbar)$ are ``generically'' isomorphic, and it is helpful to make use of the interplay between the two.
 The discussion above shows, for example, that while $\cL(\nu,\hbar)$ carries a natural vertex bi-algebra structure, $\cW(\lambda,c)$ does not, although a comultiplication can be defined upon appropriate localization.

 \subsection{ }
 \label{intro-why-chiral-group}
 In order fully to justify the name of {\em a chiral Lie group} that the title of the note seems to bestow on $\cL(\nu,\hbar)$, the latter had better be a chiral Hopf algebra, that is, one needs
 a chiral antipode $\cL(\nu,\hbar)\rightarrow\cL(\nu,\hbar)$. Note that the definition of the opposite vertex algebra $\cL(\nu,\hbar)^-$ is quite clear in the Beilinson-Drinfeld language, where
 a vertex algebra is simply a Lie algebra object in a certain category. The antipode may well exist in our situation, but we could not prove this. On the other hand it definitely
 exists in the quasiclassical limit, see (\ref{coiss-alg-inv}), which partially justifies the name. It is curious to note that the above mentioned \cite{BD}, 3.4.16 omits any mention
  of the antipode.
 
 \subsection{ }
 \label{intro-origins}
 
 There is a whole menagerie of algebras closely related to each other and to the instanton moduli problem: the deformed $W_{1+\infty}$, Cherednik's double affine
 Hecke algebra, spherical elliptic Hall algebra, the Yangian of $\widehat{gl}_1$, and whatnot, see e.g. \cite{ASch,SV}. Algebras of the type $\cW(\lambda,c)$ or
 $\cL(\nu,\hbar)$ are
 also on the list, see e.g. \cite{Pro} and references therein. More than anything else, the purpose of this note is  better to understand this emerging  picture. For example,
 the algebra of differential operators on $\BC^*$, that appeared above under the moniker of $W_{1+\infty}$ is precisely half of the Manin triple that underlines
 the Poisson-Lie structure on $G_{\Psi DO}$. To show another example, our results relate to the Yangian approach 
 to the instanton moduli \cite{MO}  as follows: according to {\em loc. cit.} a Yangian is among other things an algebra that operates on a given list
 of modules and their tensor products, and this is exactly the property  $\cL(\nu,\hbar)$  enjoys: it projects onto various conventional $W_n(\hbar)$ and, therefore,
 operates on various Fock spaces and their tensor products, where it acts via the ``comultiplication'' of Theorem~\ref{main-res}(ii). This structure is closely related
 to the way the AGT conjecture was proved in \cite{SV}.
 
 We would also like to mention the remarkably written introduction to \cite{FJMM} where the idea of analytic continuation in relation to Yangians and moduli of instantons
 is put forward.
 
 \subsection{ }
 \label{intro-thanks} 
 F.M. would like to thank Y.Soibelman for introducing him to the various aspects of the instanton moduli theory and B.Khesin for a memorable collaboration.
 
 \section{ proof of theorem~\ref{main-res}}
 \label{proof-main-res}
 \subsection{ } 
 \label{constr-vert-alge-strre-detail} Let us begin by constructing the vertex algebra structure on $\cL(\nu,\hbar)$ that appeared in sect.~\ref{intro-4}. By definition, sect.~\ref{intro-3},
 \begin{equation}
 \label{form-nuj-more-detail}
 ^nU_j(z)=\sum _{i_1<i_2<\cdots<i_a, \sum_im_i=j-a}c_{\vec{i}}^{\vec{m}}I_{i_1}(z)^{(m_1)}I_{i_2}(z)^{(m_2)}\cdots I_{i_a}(z)^{(m_a)}
 \end{equation}
 with scalar coefficients
 
\[
 c_{\vec{i}}^{\vec{m}}={i_1-1\choose m_1}{i_2-2-m_1\choose m_2}{i_3-3-m_1-m_2\choose m_3}\cdots{i_a-a-m_1-\cdots-m_{a-1}\choose m_a}.
 \]
 For example,
 \begin{eqnarray}
 ^nU_1(z)&=&\sum_jI_j(z)\nonumber\\
 ^nU_2(z)&=&\sum_{i_1<i_2}I_{i_1}(z)I_{i_2}(z)+\sum_i(i-1)I_i(z)^{(1)}\nonumber\\
 ^nU_3(z)&=&\sum_{i_1<i_2<i_3}I_{i_1}(z)I_{i_2}(z)I_{i_3}(z)+\sum_{i_1<i_2}(i_1-1)I_{i_1}(z)^{(1)}I_{i_2}(z)+\nonumber\\
 & &\sum_{i_1<i_2}(i_2-2)I_{i_1}(z)I_{i_2}(z)^{(1)}
 +\sum_i{i-1\choose 2}I_i(z)^{(2)}.\nonumber
 \end{eqnarray}
 What is important is that $c_{\vec{i}}^{\vec{m}}$ is a polynomial in the $i$-s and is independent of $n$ provided $n\geq j$ (otherwise the field $^nU_j(z)$ makes no sense.)
 
 Denote by $I_{\vec{i}}(z)^{(\vec{m})}$ the monomial that appears in the r.h.s. of (\ref{form-nuj-more-detail}).  
 
 The OPE of $^nU_i(z)\cdot ^nU_j(w)$ is a linear combination
 of the normally ordered products $: I_{\vec{p}}(z)^{(\vec{m})}I_{\vec{q}}(w)^{(\vec{l})}:$ with coefficients in $\BC[(z-w)^{-1}]$.  These coefficients are computed using the
 Wick theorem as follows: for any two pairs $(\vec{k'},\vec{r})$ and $(\vec{k''},\vec{r})$ such that the monomials $I_{\vec{p}\cup\vec{r}}(z)^{(\vec{m}\cup\vec{k'})}$ and
 $I_{\vec{q}\cup\vec{r}}(w)^{(\vec{l}\cup\vec{k''})}$
  enter expansion (\ref{form-nuj-more-detail}) of
 $^nU_i(z)$ and $^nU_j(w)$ resp. there arises a contribution equal to
 \[
c_{\vec{p}\cup\vec{r}}^{(\vec{m}\cup\vec{k'})}c_{\vec{q}\cup\vec{r}}^{(\vec{l}\cup\vec{k''})}\mbox{contr} (I_{\vec{r}}(z)^{(\vec{k'})}I_{\vec{r}}(w)^{(\vec{k''})}): I_{\vec{p}}(z)^{(\vec{m})}I_{\vec{q}}(w)^{(\vec{l})}:
 \]
 with Wick's contraction defined by iterating and differentiating (\ref{def-of-heis}):
 \[
 \mbox{contr} (I_{\vec{r}}(z)^{(\vec{k'})}I_{\vec{r}}(w)^{(\vec{k''})})=\prod_\alpha\partial_z^{k'_\alpha}\partial_w^{k''_\alpha}\frac{\hbar}{(z-w)^2}
 \]
The latter expression is independent of $\vec{r}$ and so the actual contribution involves sums of the type
\[
\sum_{\vec{r}}c_{\vec{p}\cup\vec{r}}^{(\vec{m}\cup\vec{k'})}c_{\vec{q}\cup\vec{r}}^{(\vec{l}\cup\vec{k''})},
\]
where, recall, summands are polynomials in subindices. The summation extended over those $\vec{r}=(r_1,r_2,...)$ where $r_1<r_2<\cdots n$ and where the $r_\bullet$'s are not allowed to equal
any of the entries of either $\vec{p}$ or $\vec{q}$. A moment's thought will convince the reader that such sums are polynomials in $n$ if $n$ is greater than $i+j$. (Indeed, that $\sum_{i=1}^n f(i)$ is a polynomial, provided $f$ is,
is well known, from which an inductive argument will derive the polynomiality of $\sum_{i_1<i_2<\cdots<i_k\leq n}f(i_1,...,i_k)$. The requirement that the $i$'s must miss some fixed integers is also easy
to incorporate.)

What all of this means is that the OPE of  $^nU_i(z)\cdot ^nU_j(w)$, as an element of $\cH(\hbar)((z-w))$, depends on $n$ polynomially, which is not quite what we want. What we need is to show that the
coefficients of the OPE of  $^nU_i(z)\cdot ^nU_j(w)$ w.r.t.   a basis  of $\cW_n(\hbar)$ are polynomials in $n$. In order to see this choose a basis of $\cH(\hbar)$ to consist of monomials 
$I^{(\vec{m})}_{\vec{i}}$ that correspond to the fields $I_{\vec{i}}(z)^{(\vec{m})}$ introduced above.  A basis of $\cW_n(\hbar)$ is similarly chosen to consist of the monomials determined by the fields
$^nU_i(z)$.  It is clear that the construction of $\lim{\cW_n(\hbar)}$ extends to an obvious definition of $\lim{\cH(\hbar)}$ and the basis choices agree with the limits. The desired coefficients are a solution
$\vec{x}=(x_1,x_2,...)$ of the linear system of equations
\[
A\vec{x}=[(^nU_i)_{(s)}(^nU_j)],
\]
where $(^nU_i)_{(s)}(^nU_j)$ is the $s$-product of the indicated elements, $[(^nU_i)_{(s)}(^nU_j)]$ is the corresponding coordinate vector w.r.t. the basis in $\cH(\hbar)$, $A$ is the matrix whose
 columns are labelled by the basis elements of (an appropriate conformal weight subspace of) $\cW_n(\hbar)$; written in each column are the coordinates of the respective element w.r.t. the basis of
 $\cH(\hbar)$.  This system does not have to have a solution -- in fact, typically the number of unknowns is much less than the number of equations -- but it does: this is a highly nontrivial result of Lukyanov
 \cite{Luk}. It is clear that the columns of $A$ are linearly independent: focusing on the leading terms of the entries one recognizes in them the basis of space of symmetric functions written in terms of
 elementary symmetrtic functions.  Therefore the solution can be found by using, say, Cramer's rule. The result is a polynomial in $n$ simply because so is $[(^nU_i)_{(s)}(^nU_j)]$, as discussed above.
 This concludes the construction of $\cL(\nu,\hbar)$.
 
 \subsection{ } 
 \label{proof-them-ref-i}{\em Proof of Theorem~\ref{main-res} (i).}  It is easy to see that the formula defining $\Delta$ is exactly what (\ref{def-of-alm-comulti}) equals if $m=\nu$. Therefore the assertion is an immediate
 consequence of  the construction of $\cL(\nu,\hbar)$ as ``analytic continuation'' of $\cW_n(\hbar)$ w.r.t. $n$.
 
 \subsection{ }{\em Proof of Theorem~\ref{main-res} (ii).} For $X\in \cL(\nu,\hbar)$ denote by $\langle X\rangle$ the projection of $X$ on the conformal weight 0 subspace. We need show
 that for any positive conformal weight elements $X$ and $Y$ and any  integer $i$ the projection $\langle X_{(i)}Y\rangle\in \nu \cL(\nu,\hbar)$. By definition, this amounts to lifting $X$ and $Y$ to
 $\cW_n(\hbar)$ and showing that $\langle X_{(i)}Y\rangle$ is a polynomial divisible by $n$ for $n\gg 0$.  Wick's theorem shows that for monomials $I^{(\vec{m})}_{\vec{i}}$ and  $I^{(\vec{l})}_{\vec{j}}$
 the projection $\langle( I^{(\vec{m})}_{\vec{i}})_{(k)}(I^{(\vec{l})}_{\vec{j}})\rangle$ may be nonzero only if $\vec{i}=\vec{j}$ in which case it comes from a ``total contraction,'' where each generator
 $I^{(\bullet)}_\bullet$ that enters $I^{(\vec{m})}_{\vec{i}}$ pairs with one generator that enters $I^{(\vec{l})}_{\vec{j}}$. What this means, cf. sect.~\ref{constr-vert-alge-strre-detail}, is that each such projection
 equals a sum of the type
 \[
 \sum_{1\leq i_1< i_2<\cdots<i_k\leq n}f(i_1,i_2,...,i_k)
 \]
 for some polynomial $f$. The result is of course a polynomial in $n$, but more is true:
 \[
 \sum_{1\leq i_1< i_2<\cdots<i_k\leq n}f(i_1,i_2,...,i_k)=0\mbox{ if }n=0,1,2,...,k-1.
 \]
 This elementary observation must be well known but, as the saying goes, we failed to find an appropriate reference. A simple example is, of course, this:
 \[
 \sum_{1\leq i_1< i_2<\cdots<i_k\leq n}1={n\choose k}.
 \]
 Another classic example involves the sums of fixed powers: let us denote by $\Sigma^{(k)}=\sum_{i=1}^ni^k$, and of course one has
 \[
 \Sigma^{(0)}=n,\; \Sigma^{(1)}={n+1\choose 2},
 \]
 both divisible by $n$. A standard trick
 \[
 n^k=\sum_{i=1}^n(i^k-(i-1)^k)=k\Sigma^{(k-1)}-{k\choose 2}\Sigma^{(k-2)}+{k\choose 3}\Sigma^{(k-3)}-\cdots
 \]
 allows to show by induction that each $\Sigma^{(k)}$ is a polynomial divisible by $n$. For a repeated summation one can write
 \[
 \sum_{1\leq i_1< i_2<\cdots<i_k\leq n}i_1^{a_1}i_2^{a_2}\cdots i_k^{a_k}=\sum_{i_k=k}^ni_k^{a_k}  \sum_{1\leq i_1< i_2<\cdots<i_{k-1}\leq i_k-1}i_1^{a_1}i_2^{a_2}\cdots i_{k-1}^{a_{k-1}}
 \]
 By the inductive assumption, the inner summation is a polynomial in $i_k$ vanishing at $i_k-1=0,1,2,...,k-2$, which implies that the outer summation can be extended to
 $\sum_{i_k=1}^n$ and is, by the induction basis,  a polynomial vanishing at $n=0,1,2,...,k-1$.  This completes the proof of theorem~\ref{main-res} (ii).
 
 \subsection{ }{\em Proof of Theorem~\ref{main-res} (iii).} 
 \label{proof-bial-and-quasi-limi}
 The bi-algebra assertions are quite obvious and will not be discussed. The fact that this bi-algebra structure is a quantization of the coisson structure on $\BC[G_{\Psi DO}]^c$ is also easy but merits
 a word or two. Recall that if a vector space carries vertex algebra structure over $\BC[\hbar]$ s.t. modulo $(\hbar)$ the vertex algebra becomes commutative, this quotient is automatically
 a coisson algebra w.r.t. the bracket $\{.,.\}$ defined to be $d/d\hbar[.,.]\mbox{ mod }(\hbar)$, where $[.,.]$ is the chiral bracket. The original vertex algebra is then called a quantization of this quotient.
 The relevant example is the Heisenberg algebra $\cH(\hbar)$ we started with,
 see (\ref{def-of-heis}).  Lukyanov's $\cW_n(\hbar)/(\hbar)$ is a coisson subalgebra of $\cH(\hbar)/(\hbar)$, and it is, by definition, the Miura transform of the Gelfand-Dickey $\BC[\mbox{DO}_n]^c$.
 Hence, $\cW_n(\hbar)$ is a quantization of $\BC[\mbox{DO}_n]^c$.  The reason that $\cL(\nu,\hbar)$ is a quantization of the whole $\BC[G_{\Psi DO}]^c$ is that, informally speaking,
 $\cL(\nu,\hbar)$ is an analytic continuation of the former, $\cW_n(\hbar)$, while $\BC[G_{\Psi DO}]^c$ is an analytic continuation of the latter, $\BC[\mbox{DO}_n]^c$.  Formally,
 identifying  $U_i\in \cL(\nu,\hbar)$ with its namesake $U_i\in\BC[G_{\Psi DO}]^c$, see sect.~\ref{intro-1}, one sees that by Lukyanov's result $d/d\hbar[U_i,U_j]\mbox{ mod }(\hbar)$ is precisely the
 Gelfand-Dickey bracket for $\lambda=N\in\BZ_+$ and $N$ large enough, say, $N>i+j$. On the other hand, hand both $d/d\hbar[U_i,U_j]\mbox{ mod }(\hbar)$ and Gelfand-Dickey's bracket are polynomials in $\lambda$. Therefore they are equal for all values of $\lambda$. $\qed$

 \subsection{ }{\em Proof of Theorem~\ref{main-res} (iv).}  
 
 \subsubsection{ }
 \label{disc-andys-paper}This easy-to-anticipate result is nontrivial and we heavily rely on \cite{L}. That paper deals with classification of vertex algebras strongly
  generated by a collection of fields with prescribed conformal weights. In order to state the main result needed for our purposes  consider the full subcategory $\cC$ of conformally graded vertex algebras  that satisfy the following conditions:
  
  $\bullet$ each $V\in\cC$ is strongly generated, over $\BC$, by conformal fields $u_j$  of conformal weight $j$, $j\geq2$;
  
  $\bullet$ $u_2$ is a Virasoro element with nonzero central charge $2u_{2(3)}u_2$;

  $\bullet$ $u_3$ is primary;
  
  $\bullet$ each $V\in\cC$ carries an automorphism s.t. $u_2\mapsto u_2$, $u_3\mapsto -u_3$;
  
  $\bullet$ the elements $u_2$ and $u_3$ satisfies the following non-degeneracy conditions: 
  
  \begin{equation}
  \label{non-deg-of-vira}
  u_{2(3)}u_2\neq0,
  \end{equation}

  \begin{equation}
  \label{1st-nonde}
 u_{3(5)}u_3\neq0,
   \end{equation}
   
   \begin{center} and \end{center}
  
  \begin{eqnarray}
  \label{2nd-nondege}
  u_{3(1)}u_3&=&\alpha_4 u_4+F_4(u_1,u_2,u_3),\; \alpha_4\neq0\nonumber\\
  u_{3(1)}(u_{3(1)}u_3)&=&\alpha_5 u_5+F_4(u_1,u_2,u_3,u_4),\;\alpha_5\neq0\nonumber\\
  \ldots & &\ldots\\
  \underbrace{u_{3(1)}(u_{3(1)}(\cdots(u_{3(1)}u_3)\cdots)}_{j-2}&=&\alpha_j u_j+F_j(u_1,u_2,\ldots u_{j-1}),\; \alpha_j\in\BC\setminus\{0\},\; j\geq 4,\nonumber
  \end{eqnarray}
  where $F_j(u_1,u_2,\ldots u_{j-1})$ is a polynomial expression involving indicated elements and their derivatives.
  
  The meaning of the last condition is that $u_2$ and $\tilde{u}_j$ recurrently defined by $\tilde{u}_{j+1}=u_{3(1)}\tilde{u}_j$, $\tilde{u}_3=u_3$, $j=3,4,...$, can be chosen as a strongly generating set; notice that the conformal weight of $\tilde{u}_j$ is automatically $j$.
  
  Constructed in \cite{L} is a universal vertex $\BC[\lambda,c]$-algebra $\cW(\lambda,c)$ that satisfies the following condition: for each $V\in\cC$ there is a 1-dimensional $\BC[\lambda,c]$-module $\BC_V$ and a vertex algebra isomorphism 
  \begin{equation}
  \label{andys res-1st take}
  \BC_V\otimes_{\BC[\lambda,c]}\cW(\lambda,c)\rightarrow V.
  \end{equation} 
  
  This isomorphism is not unique, but almost so: there are exactly two such isomorphisms, one being obtained from another as the composition with the above stipulated automorphism that sends
  $u_2\mapsto u_2$, $u_3\mapsto -u_3$. 
  
  Part of the definition of $\BC_V$ is easy to reproduce here: the morphism $\BC[\lambda,c]\rightarrow\BC_V$ is determined by the images of $\lambda$ and $c$; the latter has the meaning of the
  Virasoro central charge and therefore
  \begin{equation}
  \label{where-c-goe}
  c\mapsto 2u_{2(3)}u_2,
  \end{equation}
  $u_2$ being the Virasoro element of $V$. To find out where $\lambda$ is mapped the reader should consult \cite{L}, (5.6-7).
  
  \subsubsection{ }
  \label{how-bout-over-ring?} Now let $V$ be a vertex $A$-algebra s.t. $A$ is an integral domain placed in conformal weight 0. We require that 
  
  its localization to the field of quotients $Q(A)\otimes_AV$ be strongly generated by $u_j$ of conformal weight $j$,
  $j\geq 2$,
  
   the condition (\ref{2nd-nondege}) hold if we allow $\alpha_j\in A\setminus\{0\}$,
  
  the remaining conditions imposed above on $\cC$ remain true as stated. 
  
  In this case \cite{L}, cf. (\ref{andys res-1st take}), gives a ring extension $A\rightarrow\tilde{A}$, a ring morphism $\BC[\lambda,c]\rightarrow\tilde{A}$,  and a vertex algebra morphism
  \begin{equation}
  \label{andys res-2nd take}
  \tilde{A}\otimes_{\BC[\lambda,c]}\cW(\lambda,c)\longrightarrow \tilde{A}\otimes_AV.
  \end{equation}
  Here
  \[
  \tilde{A}=A[\alpha_j^{-1},(u_{2(3)}u_2)^{-1},(u_{3(5)}u_3)^{-1},j\geq4][x]/(x^2-2u_{2(3)}u_2/3u_{3(5)}u_3),
  \]
  \[
  A\longrightarrow\tilde{A}
  \]
  is the canonical map ( standard \'etale, therefore,); finally
  \[
  \BC[\lambda,c]\longrightarrow\tilde{A}\mbox{ is such that }c\mapsto 2u_{2(3)}u_2,
  \]
  cf. (\ref{where-c-goe}), and we will not describe the image of $\lambda$ again making reference to (5.6-7) in \cite{L}. The discerning reader will notice that the localization at $(u_{2(3)}u_2)^{-1},(u_{3(5)}u_3)^{-1}$ and $\alpha_j^{-1}$ has its origin
  in the non-degeneracy conditions (\ref{non-deg-of-vira},\ref{1st-nonde},\ref{2nd-nondege}), and the appearance of the square root $\sqrt{2u_{2(3)}u_2/3u_{3(5)}u_3}$ is the manifestation of there being two distinct
  morphisms $\cW(\lambda,c)\longrightarrow \tilde{A}\otimes_AV$, see the beginning of sect.~\ref{disc-andys-paper}.

  \subsubsection{ }
  \label{appl-andy-to-l}
  We would like to use (\ref{andys res-2nd take}) in order to obtain a morphism $\cW(\lambda,c)\rightarrow\cL(\nu,\hbar)$.  
  
  To begin with,  $\cL(\nu,\hbar)$ is not of the type considered in \cite{L} as it is generated by the fields of conformal weights 1,2,3,... The conformal weight 1 generator, $U_1$, produces the
  Heisenberg algebra as 
  \[
  ^nU_1(z)\cdot^nU_1(w)= \frac{n\hbar}{(z-w)^2}+\cdots,
  \]
  thanks to (\ref{def-of-heis},\ref{form-nuj-more-detail}).
  
  What one needs, therefore, is to consider the commutant of this Heisenberg in $\cL(\nu,\hbar)$.  In what follows we always work in $\cW_n(\hbar)$ and make sure that the result stabilizes as
  $n$ becomes large. The field $^nU_2(z)$ does not commute with $^nU_1(z)$ but
  \begin{equation}
  \label{vir-fixed-n}
  L(z)\stackrel{\mbox{def}}{=}\,^nU_2(z)-\frac{n-1}{n}:^nU_1(z)^nU_1(z):-\frac{n-1}{2}\,^nU_1(z)'
  \end{equation}
  does and generates the Virasoro algebra of central charge
 \begin{equation}
 \label{central-charge-fixed-n}
  c(n,\hbar)=(n-1)(1-\frac{n(n+1)}{\hbar}).
  \end{equation}
  This means that
  \[
  L(z)\stackrel{\mbox{def}}{=}U_2(z)-\frac{\nu-1}{\nu}:U_1(z)U_1(z):-\frac{\nu-1}{2}U_1(z)'
  \]
defines a Virasoro subalgebra of $\cL(\nu,\hbar)[(\nu\hbar)^{-1}]$ with central charge
  \begin{equation}
  \label{centr-charg}
  c(\nu,\hbar)=(\nu-1)(1-\frac{\nu(\nu+1)}{\hbar}).
  \end{equation}  
  
  At this point it becomes clear that
  
  (1) the localization of the commutatnt $\BC(\nu,\hbar)\otimes_{\BC[\nu,\hbar]}\mbox{Comm}(U_1(z),\cL(\nu,\hbar))$ is a vertex algebra with Virasoro element of central charge $c_{\nu,\hbar}$;
  
  (2) $\BC(\nu,\hbar)\otimes_{\BC[\nu,\hbar]}\mbox{Comm}(U_1(z),\cL(\nu,\hbar))$ is of the type considered in \cite{L}; namely, it is strongly generated by fields $\tilde{U}_j(z)$ of conformal dimension $j$, $j\geq 2$, and, in fact, one can choose $\tilde{U}_j=U_j+F_j(U_1,...,U_{j-1})$;
  
  (3)  all the conditions imposed on $V$ in sect.~\ref{how-bout-over-ring?} in order that (\ref{andys res-2nd take}) exist are satisfied.
  
  Item (3), the only one worth being commented on, is valid because the conditions imposed on the category $\CC$ at the beginning of sect.~\ref{disc-andys-paper}
  are axiomatizations of the known properties of $\cW_n(\hbar)$ of which $\cL(\nu,\hbar)$
  is but an analytic continuation. For example, the existence of the automorphism $u_2\mapsto u_2,u_3\mapsto-u_3$ follows from \cite{Luk}, Proposition 3. This guarantees the existence of a morphism
  \[
\cW(\lambda,c)\rightarrow\BC(\nu,\hbar)\otimes_{\BC[\nu,\hbar]}\cL(\nu,\hbar),
\]
We can, however, do much better and, in particular, avoid the passage to the total field of fractions. We are about to see that for our purposes it suffices to compute the next generating field,
$\tilde{U}_3$.

A bit of a computational  effort will show that the field
\begin{eqnarray}
\label{tilde-u-3}
^n\tilde{U}_3(z)&\stackrel{\mbox{def}}{=} &^nU_3(z)-\frac{n-2}{2}\,^nU_2(z)'+\frac{(n-1)(n-2)}{12}\,^nU_1(z)^{(2)}\nonumber\\
&+&\frac{(n-1)(n-2)}{2n}:^nU_1(z)'\,^nU_1(z):
-\frac{n-2}{n}:^nU_1(z)\,^nU_2(z):\nonumber\\
&+&\frac{(n-1)(n-2)}{3n^2}:^nU_1(z)^nU_1(z)^nU_1(z):\nonumber
\end{eqnarray}
commutes with $^{n}U_1(z)$ and is primary with respect to $L(z)$.  This allows to define $\tilde{U}_3\in\cL(\nu,\hbar)[\nu^{-1}]$ by the familiar device of replacing $n$ with $\nu$ in the formula above.

What remains is to verify is  that the recurrently defined $\tilde{U}_3$, $\tilde{U}_4=\tilde{U}_{3(1)}\tilde{U}_3$, $\tilde{U}_{r+1}=\tilde{U}_{3(1)}\tilde{U}_r$, $r\geq 3$ strongly generate 
$\cL(\nu,\hbar)$, upon appropriate localization.  This computation is a pleasing exercise on arguments of sect.~\ref{constr-vert-alge-strre-detail} . 

By construction, the indicated elements are uniquely determined by their leading terms, which are symmetric polynomials in the Heisenberg algebra generators $I_1,I_2,...$ with no derivatives involved.
It is clear that the leading term of any element of $\cL(\nu,\hbar)$ is a symmetric polynomial in the indicated variables, and, by definition, the leading term of the original $U_r$ is (identified with) the
$r$-th {\em elementary} symmetric function $e_r$. What we need is to show that the leading term of $\tilde{U}_r$ has the form $\alpha_re_r+F_r(e_1,e_2,...,e_{r-1})$ with $\alpha_r\neq0$.
Since we are only interested in the leading term of $\tilde{U}_{3(1)}\tilde{U}_r$, it suffices to compute $e_{3(1)}e_r$. A straightforward application of Wick's theorem gives
\[
e_{3(1)}e_r={r+1\choose 2}(n-r+1)e_{r+1}+(r-1)(n-r)m_{(2,1^{r-1})}+(n-r+1)m_{(2,2,1^{r-3})},
\]
where we use some standard notation from theory of symmetric functions, see \cite{Mac}. Re-writing this expression as a polynomial in elementary symmetric functions one obtains
\[
e_{3(1)}e_r=-(r+1) e_{r+1}+\cdots
\]
and so
\[
\tilde{U}_r=\frac{(-1)^{r+1}r!}{6}e_r+\cdots
\]
This means that, somewhat surprisingly, no localization is needed to ensure that the fields $\tilde{U}_2\stackrel{\mbox{def}}{=}L(z)$ and $\tilde{U}_3$ generate the  algebra.

To conclude:

\begin{lem}
\label{form-ands-morphism-lemma}
There is an injection
\[
\cW(\lambda,c)\longrightarrow\cL(\nu,\hbar)[S][\sqrt{\frac{2\nu\hbar}{(\nu-2)(4\hbar-2\nu-2\nu^2)}}],
\]
where $S$ is the multiplicative set generated by $\{\hbar^{-1},(\nu-1)^{-1},(\nu-2)^{-1},(\hbar-\nu-\nu^2)^{-1},(4\hbar-2\nu-\nu^2)^{-1}\}$.
\end{lem}
{\em Proof.} The existence of the morphism is clear from the argument above, the injectivity follows from the -- generic in $\lambda,c$ -- simplicity of $\cW(\lambda,c)$, and what needs to be explained is the exact form of localization that is asserted. These formulas are
not very important for us here, but we will sketch the proof nevertheless. Denote the generating fields of $\cW(\lambda,c)$ by $W_j(z)$, $j=2,3,4,...$. The morphism is determined
by the assignment
\[
W_2\mapsto x \tilde{U}_2,\; W_3\mapsto y\tilde{U}_3\; c\mapsto c(\nu,\hbar),\lambda\mapsto\lambda(\nu,\hbar),
\]
and what we need to determine is 4 elements, $x,y,c(\nu,\hbar),\lambda(\nu,\hbar)$, of an appropriate (also to be determined) extension of $\BC[\nu,\hbar]$.

By definition, $W_2$ is a Virasoro element with central charge $c$ so that $W_{2(3)}W_2=c/2$ A glance at (\ref{centr-charg}) shows that 
\[
c\mapsto c(\nu,\hbar)=(\nu-1)(1-\frac{\nu(\nu+1)}{\hbar})=\frac{(\nu-1)(\hbar-\nu-\nu^2)}{\hbar}),
\]
which recovers one of our formulas.

 Being a conformal weight 3 primary field, $W_3$ is unique up to proportionality. The normalization chosen in \cite{L} requires that
\[
W_{3(5)}W_3=\frac{c}{3}.
\]
On the other hand $\tilde{U}_{3}$, (\ref{tilde-u-3}), satisfies
\[
\tilde{U}_{3(5)}\tilde{U}_{3}=\frac{(\nu-1)(\nu-2)(\hbar-\nu-\nu^2)(4\hbar-2\nu-\nu^2)}{6\nu\hbar^2},
\]
cf. \cite{Pro}, (169).

This forces the assignment
\[
W_3\mapsto\sqrt{\frac{2\nu\hbar}{(\nu-2)(4\hbar-2\nu-\nu^2)}}\tilde{U}_3.
\]
It remains to figure out the image of $\lambda$. The easiest -- and the most useful one, as we shall soon see -- way to do so is to cite another result of \cite{L}. Theorem~7.4 of {\em loc. cit.} claims that
the ordinary W-algebra is a homomorphic image of $\cW(\lambda,c)$. Namely, if we introduce localization
$\cW(\lambda,c)[t^{-1}]$, $t=(-2+2c-n+cn+3n^2)$, then we obtain  a surjection
\begin{equation}
\label{morph-onto-ordina-w-from- andys}
\cW(\lambda,c)[t^{-1}]\twoheadrightarrow \cW(sl_n)
\end{equation}
so that
\begin{equation}
\label{morph-onto-ordina-w-from- andys-lambda}
c\mapsto c,\;\lambda\mapsto \frac{(n-1)(n+1)}{(n-2)(-2+2c-n+cn+3n^2)},
\end{equation}
where we have taken the liberty of identifying the central charge of the canonical Virasoro element in   $\cW(\lambda,c)$ and $\cW(sl_n)$ with $c$.
Replacing in the last expression $n$ with $\nu$ and $c$ with the above computed  $c(\nu,\hbar)$ implies that in our case
\begin{equation}
\label{formu-for-lambda-nu-h}
\lambda\mapsto\lambda(\nu,\hbar)\stackrel{\mbox{def}}{=}\frac{\hbar}{(\nu-2)(4\hbar-2\nu-\nu^2)}.
\end{equation}
$\qed$

{\em Remark.} For the sake of the demanding reader here is a quick reminder on $\cW(sl_n)$. It is customary to introduce the W-algebra as the quantum Hamiltonian reduction of the corresponding
affine Lie algebra, but for our purposes  original Lukyanov's approach \cite{Luk} is much more convenient.  $\cW_n(\hbar)$ contains the Heisenberg that is generated by
the field $^nU_1(z)$ and $\cW(sl_n)$ is its commutant.  This definition brings us back to the beginning of the current section: for example, $\cW(sl_n)$ has a unique up to proportionality conformal weight 2
element, given by (\ref{vir-fixed-n}), which generates the Virasoro algebra with central charge (\ref{central-charge-fixed-n}).

\subsubsection{ }
\label{conclu-of-iii}
It is easy now  to conclude the proof of Theorem~\ref{main-res} (iii). We have 2 morphisms,
\[
\cW(\lambda,c)[t^{-1}]\longrightarrow\cL(\nu,\hbar)[S][\sqrt{\frac{2\nu\hbar}{(\nu-2)(4\hbar-2\nu-2\nu^2)}}]
\]
and, see (\ref{morph-onto-ordina-w-from- andys}),
\[
\cW(\lambda,c)[t^{-1}]/(\lambda-\frac{(n-1)(n+1)}{(n-2)(-2+2c-n+cn+3n^2)})\twoheadrightarrow \cW(sl_n).
\]
Tensoring with the 1-dimensional Heisenberg algebra, $H$, these become
\begin{equation}
\label{fro-wto-luky}
H\otimes\cW(\lambda,c)[t^{-1}]\longrightarrow\cL(\nu,\hbar)[S][\sqrt{\frac{2\nu\hbar}{(\nu-2)(4\hbar-2\nu-2\nu^2)}}]
\end{equation}
and
\begin{equation}
\label{fromm-w-to-w}
H\otimes\cW(\lambda,c)[t^{-1}]/(\lambda-\frac{(n-1)(n+1)}{(n-2)(-2+2c-n+cn+3n^2)})\twoheadrightarrow \cW_n(\hbar).
\end{equation}

By construction, morphism (\ref{fro-wto-luky}) sends a basis of its source to a basis of its target; it would be  an isomorphism were the scalars of the former the same as the scalars of the latter.
They are not and to invert (\ref{fro-wto-luky}) one needs to extend scalars on the left as follows. Let
\[
A=\BC[\nu,\hbar][S][\sqrt{\frac{2\nu\hbar}{(\nu-2)(4\hbar-2\nu-2\nu^2)}}].
\]
This ring is a receptacle of  $\BC[\lambda,c]$ under  (\ref{fro-wto-luky}); therefore  we obtain the desired vertex algebra morphism
\begin{equation}
\label{fro-wto-luky-inversi}
A\otimes_{\BC[\lambda,c,t^{-1}]}H\otimes\cW(\lambda,c)[t^{-1}]\longleftarrow\cL(\nu,\hbar).
\end{equation}

Upon the performed ring extension map (\ref{fromm-w-to-w}) becomes
\begin{equation}
\label{fromm-w-to-w-2}
A\otimes_{\BC[\lambda,c,t^{-1}]}H\otimes\cW(\lambda,c)[t^{-1}]/(\lambda-\frac{(n-1)(n+1)}{(n-2)(-2+2c-n+cn+3n^2)})\twoheadrightarrow \cW_n(\hbar).
\end{equation}

If we replace $\lambda$ with $\lambda(\nu,\hbar)$, see (\ref{formu-for-lambda-nu-h}), and $c$ with $c(\nu,\hbar)$, see (\ref{centr-charg}), in 
\[\lambda-\frac{(n-1)(n+1)}{(n-2)(-2+2c-n+cn+3n^2)},
\]
 the latter expression becomes
\[
-\frac{3(\nu-n)(\hbar n^2-n\nu-\nu^2)}{(n-2)(4\hbar-2\nu-\nu^2)},
\]
and so the principal ideal 
\[
(\lambda-\frac{(n-1)(n+1)}{(n-2)(-2+2c-n+cn+3n^2)})
\]
 equals the principal ideal
\[
(-\frac{3(\nu-n)(\hbar n^2-n\nu-\nu^2)}{(n-2)(4\hbar-2\nu-\nu^2)}).
\]

 Thanks to the appearance of  the factor of $(\nu-n)$ in this expression\footnote{which, by the way, is not at all surprising as the vanishing of this expression at $\nu=n$
 is built into the definition of $\lambda(\nu,\hbar)$, see (\ref{formu-for-lambda-nu-h})}, the composition of (\ref{fro-wto-luky-inversi}) and (\ref{fromm-w-to-w-2})
 \[
 \cL(\nu,\hbar)\longrightarrow \cW_n(\hbar)
 \]
 factors through the projection 
 \[
 \cL(\nu,\hbar)\twoheadrightarrow \cL(\nu,\hbar)/(\nu-n)
 \]
 and gives the desired 
 \[
 \cL(\nu,\hbar)/(\nu-n)\longrightarrow \cW_n(\hbar).\;\; \qed
 \]

\footnotesize{A.L.: Department of Mathematics, University of Denver, Denver, CO 80208, USA.
E-mail: andrew.linshaw@du.edu}

\footnotesize{F.M.: Department of Mathematics, University of Southern California, Los Angeles, CA 90089, USA.
E-mail:fmalikov@usc.edu}

 \end{document}